%% file: main.tex
\documentclass[11pt]{article}
\usepackage[english]{style}
\usepackage{latexsym,amsmath,amssymb,amsfonts}
\usepackage{amsmath,amscd}

\author{%
Damien Calaque\footnote{D.C. is on leave of absence from Institut Camille Jordan UMR5208, 
Universit\'e Lyon 1, 69622 Villeurbanne, France.}, 
Andrei C\u ald\u araru, 
Junwu Tu}
\title{$\PBW$ for an inclusion of Lie algebras}
\date{}
\input{defs}
\begin{document}
\maketitle
\input{abstract}
\input{intro}
\input{hkr}

\input{fon}
\input{pbw}
\input{twisted}
\input{example}
\input{appendix}

\input{biblio}
\begin{center}
\begin{small}
\begin{tabular}{ll}
Damien Calaque & Andrei C\u ald\u araru, Junwu Tu \\
Departement Mathematik & Mathematics Department\\
ETH Z\"urich & University of
Wisconsin--Madison\\
R\"amistrasse 101 & 480 Lincoln Drive\\
8092 Z\"urich & Madison, WI 53706\\
Switzerland & USA \\
{\tt damien.calaque@math.ethz.fr} &  
{\tt andreic@math.wisc.edu} \\
& {\tt tu@math.wisc.edu}
\end{tabular}
\end{small}
\end{center}

\end{document}

%% file: defs.tex
\DeclareFontFamily{U}{rsf}{}
\DeclareFontShape{U}{rsf}{m}{n}{
  <5> <6> rsfs5 <7> <8> <9> rsfs7 <10-> rsfs10}{}
\DeclareMathAlphabet{\mathscr}{U}{rsf}{m}{n}

\DeclareMathAlphabet{\mathgth}{U}{euf}{m}{n}

\DeclareFontFamily{U}{cyr}{}
\DeclareFontShape{U}{cyr}{m}{n}{
  <5> wncyr5 <6> wncyr6 <7> wncyr7 <8> wncyr8 <9> wncyr9 <10-> wncyr10}{}
\DeclareMathAlphabet{\mathcyr}{U}{cyr}{m}{n}

\input cyracc.def

\makeatletter
\def\operator@font{\sf}
\makeatother

\newcommand{\gh}{\mathgth{h}}
\newcommand{\gog}{\mathgth{g}}
\newcommand{\gn}{\mathgth{n}}

\newcommand{\Ug}{{\mathbf U}(\gog)}

\newcommand{\cH}{{\mathscr H}}

\newcommand{\cO}{{\mathscr O}}

\newcommand{\PBW}{{\mathsf{PBW}}}

\newcommand{\D}{{\mathbf D}}

\newcommand{\chk}{{\scriptscriptstyle\vee}}

\newcommand{\Ld}{\mathbf{L}}

\newcommand{\op}{{\mathsf{op}}}

\newcommand{\bbk}{\mathbf{k}}

\newcommand{\bone}{{\mathbf 1}}

\renewcommand{\S}{{\mathbf{S}}}

\DeclareMathOperator{\ad}{ad}

\DeclareMathOperator{\Hom}{Hom}

\DeclareMathOperator{\id}{id}

\DeclareMathOperator{\gr}{gr}

\DeclareMathOperator{\Ext}{Ext}

\DeclareMathOperator{\Res}{Res}
\DeclareMathOperator{\Ind}{Ind}

\newcommand{\bS}{\mathbf{S}}
\newcommand{\T}{\mathbf{T}}
\newcommand{\U}{\mathbf{U}}

\newcommand{\ra}{\rightarrow}

\newcommand{\gMod}{\mbox{-}\mathgth{Mod}}

\newcommand{\iso}{\cong}

\newcommand{\pj}{\mathbf{P}}
\newcommand{\img}{\mathsf{Im}~}

\renewcommand{\phi}{\varphi}

%% file: abstract.tex
\begin{abstract}
  Let $\gh \subset \gog$ be an inclusion of Lie algebras with quotient
  $\gh$-module $\gn$.  There is a natural degree filtration on the
  $\gh$-module $\Ug/\Ug\gh$ whose associated graded $\gh$-module is
  isomorphic to $\S(\gn)$.  We give a necessary and sufficient
  condition for the existence of a splitting of this filtration.  In
  turn such a splitting yields an isomorphism between the
  $\gh$-modules $\Ug/\Ug\gh$ and $\S(\gn)$.  For the diagonal
  embedding $\gh \subset \gh\oplus \gh$ the condition is automatically
  satisfied and we recover the classical Poincar\'e-Birkhoff-Witt
  theorem.

  The main theorem and its proof are direct translations of results in
  algebraic geometry, obtained using an {\it ad hoc} dictionary.  This
  suggests the existence of a unified framework allowing the
  simultaneous study of Lie algebras and of algebraic varieties, and a
  closely related work in this direction is on the way.
\end{abstract}

%% file: intro.tex
\section{Introduction}
\label{sec:intro}

\subsection{The aim}

Let $\gh\hookrightarrow \gog$ be an inclusion of Lie algebras.  Denote
by $\gn$ the quotient $\gog/\gh$.  The quotient $\U(\gog)/\U(\gog)
\gh$ of $\U(\gog)$ by the left ideal generated by $\gh$ is naturally
an $\gh$-representation.  The main purpose of this paper is to answer
the following question (the PBW problem):

\begin{center}
  \textit{When is $\U(\gog)/\U(\gog) \gh$ isomorphic to $\bS(\gn)$ as
    $\gh$-representations?}
\end{center}

\noindent
A more precise way of stating the above question is the following. The representation $\Ug/\Ug\gh$ admits a natural filtration by
$\gh$-modules whose associated graded $\gh$-module is $\S(\gn)$.  We
ask for a necessary and sufficient condition for this filtration to
split.

This question is important in deformation quantization, as the space
of $\gh$-invariants $\left ( \Ug/\Ug\gh \right )^\gh$ can be given a
natural structure of algebra by identifying it with the space of
invariant differential operators on a homogeneous space~\cite{Tor}.
An open conjecture of Duflo is concerned with understanding the center
of this algebra in terms of the Poisson center of $\S(\gn)^\gh$, which
is thought of as the algebra of functions on a Poisson manifold
obtained via reduction through the moment map $\gog^\chk \ra
\gh^\chk$.  In order for this conjecture to make sense one needs to be
in a situation where the PBW isomorphism holds.  Traditionally this is
achieved by assuming that the inclusion $\gh \hookrightarrow \gog$
splits as a map of $\gh$-modules.  We will see that this condition is
unnecessarily restrictive: there are many pairs of Lie algebras for
which there is a PBW isomorphism (and hence it makes sense to study
the Duflo problem), but for which the inclusion $\gh\hookrightarrow
\gog$ does not split.

\subsection{An analogous problem in algebraic geometry}

Kontsevich and Kapranov~\cite{Ka} had the insight that we can view the
shifted tangent sheaf $T_Y[-1]$ of a smooth algebraic variety
$Y$ as a Lie algebra object in the derived category $\D(Y)$ of
coherent sheaves on $Y$, with bracket given by the Atiyah class of the
tangent sheaf.  Moreover, the Atiyah class of any object in $\D(Y)$
gives it the structure of module over this Lie algebra object (see for
example~\cite{RW}).  Loosely speaking $\D(Y)$ can be regarded as the
category of representations of the shifted tangent sheaf.  The role of
the trivial representation is played by the structure sheaf $\cO_Y$.

An embedding $i:X\hookrightarrow Y$ of smooth algebraic varieties can
be thought of as giving rise to an inclusion of Lie algebra objects in
$\D(X)$
\[ \mathfrak h=T_X[-1]\hookrightarrow i^*T_Y[-1]=\mathfrak{g}. \] 
If $E$ is an object in $\D(Y)$ then the Atiyah class of the
restriction $i^*E$ of $E$ to $X$ is precisely the composite of the
above inclusion of Lie algebras with the restriction to $X$ of the
Atiyah class of $E$.  In other words the functor 
\[ i^*:\D(Y) \ra \D(X) \]
can be interpreted as the restriction functor 
\[ \Res:\gog\gMod \ra \gh\gMod. \] 
(We think of all our functors between derived categories as being
implicitly derived, so we write $i^*$ instead of $\Ld i^*$, etc.)

We now see a dictionary emerging between the worlds of Lie theory and
of algebraic geometry. We can use this dictionary to translate naively
the PBW question into a problem in algebraic geometry.  The following
concepts are matched by this dictionary:

\begin{center}
\begin{tabular}{l|l}
Lie theory & Algebraic geometry\\
\hline
Lie algebras $\gh$, $\gog$ & varieties $X$, $Y$, $\gh = T_X[-1]$,
$\gog = T_Y[-1]$ \\ 
inclusion $\gh\hookrightarrow \gog$ & closed embedding
$i:X\hookrightarrow Y$ \\
$\gh\gMod$, $\gog\gMod$ & $\D(X)$, $\D(Y)$ \\
$\bone_\gh\in \gh\gMod$ & $\cO_X \in \D(X)$ \\
$\Res:\gog\gMod \ra \gh\gMod$ & $i^*: \D(Y) \ra \D(X)$
\\
$\Ind:\gh\gMod \ra \gog\gMod$ & $i_!:\D(X) \ra \D(Y)$
\end{tabular}
\end{center}

\noindent
The last line is motivated by the fact that the induction functor
$\Ind$ is the left adjoint of the restriction functor, hence in the
right column we take the left adjoint $i_!$ of the pull-back functor,
which exists for a closed embedding $i$ of smooth varieties.

In representation-theoretic language the $\gh$-representation
$\Ug/\Ug\gh$ arises as
\[ \Ug/\Ug\gh = \Res \Ind \bone_\gh \in \gh\gMod. \] 
Using the dictionary the latter corresponds to the object $i^* i_!
\cO_X$ of the derived category $\D(X)$.  Any object $E$ of $\D(X)$
admits a natural filtration by successive truncations $\tau^{\geq k}
E$ whose $k$-th ``quotient'' is the cohomology sheaf $\cH^k(E)[-k]$.
An easy local calculation shows that for $E = i^* i_! \cO_X$ we have
\[ \cH^k(i^* i_! \cO_X) = \wedge^k N \]
where $N$ is the normal bundle of $X$ in $Y$.  Thus the associated
graded object of $i^*i_! \cO_X$ is precisely
\[ \gr (i^* i_! \cO_X) = \bigoplus_k \wedge^k N[-k] = \S(N[-1]). \]
Since $N[-1] = T_Y[-1]|_X/T_X[-1]$ corresponds via the dictionary to 
$n= \gog/\gh$, this is the precise analogue of the statement that
$\Ug/\Ug\gh$ admits a filtration whose associated graded is
\[ \gr \left (\Ug/\Ug\gh \right ) = \S(n). \]
 
The PBW question translates into the following question about a
closed embedding $i:X\hookrightarrow Y$.
\begin{center}
\textit{When is $i^*i_!\cO_X$ isomorphic to $\S(N[-1])$ in
  $\D(X)$?}
\end{center}
Just like in the usual PBW problem, this question is better phrased by
asking when the above filtration on $i^*i_!\cO_X$ splits.  This
question was addressed and solved recently by D.~Arinkin and the second
author in~\cite{AC}, where they prove the following result.
\medskip

\label{thm:alggeom}
\noindent \textbf{Theorem~1.2.}
{\em 
Let $X^{(1)}$ be the first infinitesimal neighborhood of $X$ in $Y$,
$X\hookrightarrow X^{(1)}$.  The following are equivalent:
\medskip
\begin{enumerate}
\item the truncation filtration on $i^*i_!\mathcal O_X$ splits, giving
  rise to an isomorphism
  \[ i^*i_! \cO_X \iso \S(N[-1]); \]
\item the class $\alpha$ is trivial, where
  \[ \alpha\in\Ext^1_X(N[-1]^{\otimes 2},N[-1]) \] 
  is obtained by composing the class of the normal bundle exact sequence
  with the Atiyah class of the normal bundle $N$;
\item the vector bundle $N[-1]$ admits an extension to $X^{(1)}$.
\end{enumerate}
}
\medskip

\noindent
It is worth noting that there are many cases where the short exact sequence
\[ 0 \ra T_X \ra T_Y|_X \ra N \ra 0\] 
does not split but the obstruction $\alpha$ is nonetheless trivial.
For example this is the case when $X$ is any non-linear hypersurface
in $Y = \pj^n$.

\subsection{The result}

Our main result is the following translation of the above theorem.
\medskip

\noindent\label{thm:mainthm}{\bf Theorem~1.3.}
{\em 
  There exists a Lie algebra $\gh^{(1)}$, containing $\gh$ as a
  Lie subalgebra, such that the following are equivalent:
\medskip
\begin{enumerate}
\item the natural filtration on $\Ug/\Ug \gh$ splits, giving rise to an isomorphism of $\gh$-modules
  \[ \Ug/\Ug\gh\iso \S(\gn); \]
\item the class $\alpha$ is trivial, where
  $\alpha\in\Ext^1_\gh(\gn^{\otimes2},\gn)$ is obtained by composing the class of
\[ 0\longrightarrow \mathfrak h\longrightarrow\mathfrak
g\longrightarrow\mathfrak n\longrightarrow 0 \]
with the $\gh$-action;
\item the $\gh$-representation $\gn$ admits an extension to $\gh^{(1)}$.
\end{enumerate}
}
\medskip

\noindent
Observe that in the algebro-geometric context $X^{(1)}$ is singular,
even though $X$ and $Y$ are smooth. It turns out that the correct
notion of tangent bundle for such singular spaces is that of {\em
  tangent complex}, see~\cite{Illusie}.  Consequently the Lie algebra
$\gh^{(1)}$ shall be defined following the analogy with the tangent
complex of $X^{(1)}$, using insight from Koszul duality.  The details
are presented in Section~\ref{sec:hkr}.

The paper is organized as follows. Section~\ref{sec:hkr} is devoted to
the definition of the ``first order neighborhood Lie algebra''
$\gh^{(1)}$ and to the proof that an $\gh$-module $E$ admits an
extension to $\gh^{(1)}$ if and only if a certain class $\alpha_E\in
Ext^1(\gn\otimes E,E)$ is trivial.  In Section 3 we prove a variant of
our main theorem for the inclusion $\gh\hookrightarrow\gh^{(1)}$.
More precisely, we prove that a natural filtration on
$U(\gh^{(1)})/U(\gh^{(1)})\gh$ splits if and only if the class
$\alpha:=\alpha_\gn$ is trivial.  The following section is concerned
with proving Theorem~\ref{thm:mainthm}.  While the main theorem is
concerned with the study of the induction-restriction of the trivial
representation, we can deduce from this case a general result for any
$\gh$-representation.  A sketch of the general case is contained in
Section~\ref{sec:twisted}.  We conclude the paper with a very short
section in which we give a simple example of a pair of Lie algebras
for which the class $\alpha$ is non-trivial.
Appendix~\ref{sec:appendix} contains a proof of the fact that a
certain algebra used in the definition of $\gh^{(1)}$ is Koszul.
\medskip

\noindent
\textbf{Assumptions.}
In what follows all algebraic structures are considered over a given
field $\bbk$.  For the main result we need to assume that the
characteristic of $\bbk$ is zero, but all other results hold without
this assumption. 
\medskip

\noindent
\textbf{Acknowledgements.} We are grateful to D.\ Grinberg for his
careful reading of the paper and for pointing out two important
mistakes in a previous version.  Discussions with D.\ Grinberg, G.\
Felder, and C.\ Rossi helped improve the exposition and clarify the
statement and proof of Theorem~\ref{thm:twisted}. We also extend
thanks to M.\ Duflo and C.\ Rossi who commented on early versions of
the result and provided the encouragement to write this paper.  The
second and third authors were partially supported by the National
Science Foundation under Grant No. DMS-0901224.

%% file: hkr.tex
\section{First order neighborhood Lie algebras}
\label{sec:hkr} 

Let $\gh\hookrightarrow \gog$ be an inclusion of Lie algebras and
denote by $\gn$ the quotient $\gh$-module $\gog/\gh$.  In this section
we define the obstruction class $\alpha$ and the first order neighborhood Lie algebra $\gh^{(1)}$ that appear in Theorem~\ref{thm:mainthm}.  Then we prove that the class $\alpha$ is trivial if and only if $\gn$ admits an extension to $\gh^{(1)}$.

\paragraph 
\label{subsec:alpha}
\textbf{The extension class $\alpha$.} We begin with the definition of
the extension class $\alpha$ that appears in the statement of
Theorem~\ref{thm:mainthm}.  Consider the short exact sequence of
$\gh$-modules
\begin{equation}
\label{ses:1} 0\rightarrow \gh \rightarrow \gog \rightarrow \gn
\rightarrow 0.
\end{equation} 
Let $E$ be an $\gh$-module.  Tensoring \eqref{ses:1} with $E$ yields
the sequence
\begin{equation}
\label{ses:2} 0\rightarrow \gh\otimes E \rightarrow \gog \otimes E
\rightarrow \gn\otimes E \rightarrow 0
\end{equation} 
which remains exact because the tensor product of representations is
the tensor product of vector spaces endowed with the $\gh$-module
structure given by the Leibniz rule.  The extension class
of~(\ref{ses:2}) is a map $\gn\otimes E \rightarrow \gh\otimes E[1]$
in the derived category of $\gh$-representations, which can be
post-composed with the action map $\gh\otimes E\rightarrow E$ to give
the map
\[ \alpha_E:\gn\otimes E \rightarrow E[1]. \] 
Equivalently, we can define $\alpha_E$ as the class in
$\Ext^1_\gh(\gn\otimes E, E)$ corresponding to the bottom extension in
the diagram below:
\begin{equation}
\label{ses:3}
\begin{CD} 0 @>>> \gh\otimes E @>>> \gog\otimes E @>>> \gn\otimes E
@>>> 0 \\ @.  @VVV @VVV @| @. \\ 0 @>>> E @>>> Q @>>> \gn\otimes E
@>>> 0.
\end{CD}
\end{equation} 
Here the $\gh$-module $Q$ is obtained by push-out in the first square
of the above diagram.  Explicitly, it is given by
\[ Q= E\oplus (\gog\otimes E) / \langle (h(x),0)-(0,h\otimes
x)\rangle \] 
where for $h\in\gh$ and $x\in E$ we have denoted by $h(x)$ the action
of $h$ on $x$ and $h\otimes x$ is viewed as an element of $\gog\otimes
E$ via the inclusion of $\gh$ into $\gog$.

We will be particularly interested in the class $\alpha_\gn$ of
the $\gh$-module $\gn$.  This special class will be denoted simply by $\alpha$.

\paragraph 
\textbf{The first order neighborhood Lie algebra $\gh^{(1)}$.}
Consider the Lie algebra $\gh^{(1)}$ defined by
\[ \gh^{(1)} := \Ld(\gog)/\langle [h,g]-[h,g]_\gog~|~h\in\gh, g\in
\gog \rangle \] 
where $\Ld(\gog)$ denotes the free Lie algebra generated by the vector
space $\gog$ and $\langle\rangle$ stands for ``Lie ideal generated by''. 
More precisely $\gh^{(1)}$ is the quotient of $\Ld(\gog)$ in
which the bracket between elements of $\gh$ and $\gog$ has been
identified with the original one in $\gog$. Note that to define the
Lie algebra $\gh^{(1)}$ we do not need $\gog$ to be a Lie algebra. The
precise weaker condition for which this construction makes sense is given in
Lemma~\ref{lem:obs} below.

There are natural maps of Lie algebras
\[ \gh \hookrightarrow \gh^{(1)} \mbox{ and } \gh^{(1)}\rightarrow
\gog\] 
which factor the original inclusion $\gh\hookrightarrow \gog$.  Given
an $\gh$-representation $E$ we can ask whether $E$ extends to a
representation of $\gh^{(1)}$. In other words we ask if on the vector
space $E$ we can find an $\gh^{(1)}$-module structure whose
restriction to $\gh$ via the map $\gh\ra \gh^{(1)}$ is the original
one.  The following lemma shows that this is the case if and only if
$\alpha_E=0$. We state the lemma in a slightly greater generality.

\begin{Lemma}
\label{lem:obs} 
Let $\gh$ be a Lie algebra and let $\gog$ be an $\gh$-module that
contains $\gh$ as an $\gh$-submodule. An $\gh$-module $E$ is the
restriction of an $\gh^{(1)}$-module if and only if its class
$\alpha_E$ is trivial.
\end{Lemma}

\begin{Proof} 
We begin with the if part. Assume that the class $\alpha_E$ is
trivial. This implies that the sequence \eqref{ses:3} splits in the
category of $\gh$-modules. Thus we get a map $j:Q\ra E$ of
$\gh$-modules that splits the canonical map $E\ra Q$. Pre-composing
$j$ with the middle vertical map in \eqref{ses:3} yields a map of
$\gh$-modules
\[ \rho: \gog\otimes E \rightarrow E.\] 
This map does not define a representation of $\gog$ on $E$, but it
certainly defines a representation of $\Ld(\gog)$ by the universal
property of $\Ld(\gog)$. The fact that $\rho$ respects the $\gh$
structure translates into the fact that $\langle
[h,g]-[h,g]_\gog\rangle $ is in the kernel of this
representation. Thus $\rho$ gives an $\gh^{(1)}$-module structure on
$E$ which lifts the original $\gh$-module structure because the first
square in \eqref{ses:3} commutes.

For the only if part assume we have an $\gh^{(1)}$-module structure on
$E$ that lifts the $\gh$ structure. Again denote this
action by $\rho$. We can use the explicit description of $Q$
above to define a splitting
\[ (x,g\otimes y)\mapsto (x+\rho(g)(y)) .\] 
This map is obviously a splitting and it respects the $\gh$-module
structure because $\langle [h,g]-[h,g]_\gog\rangle$ is in the kernel
of the representation $\rho$.  \qed
\end{Proof}

%% file: fon.tex
\section{PBW for inclusions into first order neighborhoods}
\label{sec:fon}

In this section we study the $\PBW$ property for the inclusion
$j:\gh\hookrightarrow\gh^{(1)}$ of $\gh$ into its first order
neighborhood Lie algebra $\gh^{(1)}$. We prove that the $\PBW$ theorem
holds if and only if the extension class $\alpha$, defined in
Section~\ref{sec:hkr}, vanishes.

\paragraph
We begin with some notation that will be used. Denote the Lie algebra
inclusion $\gh\hookrightarrow \gog$ by $i$. Denote the natural maps of
Lie algebras $\gh\rightarrow \gh^{(1)}$ and $\gh^{(1)}\rightarrow
\gog$ by $j$ and $k$ respectively so that $i=k\circ j$.  Denote by
$i^*$ the restriction functor from $\gog$-modules to $\gh$-modules and
by $i_!$ the induction functor in the reverse direction. Thus we have
the adjunction $i_!\dashv i^*$. We also have similar functors and
adjunctions for the maps $j$ and $k$.  Finally we denote the
$1$-dimensional trivial representation of the Lie algebra $\gh$ by
$\bone_\gh$.

\paragraph
The goal of the current paper is to understand $\PBW$ properties for
\[
i^*i_!(\bone_\gh)=\U(\gog)\otimes_{\U(\gh)}\bone_\gh=\U(\gog)/\U(\gog)\gh. \]
In this section we study the object $j^*j_!(\bone_\gh)$ which is
easier to understand. This representation can be described as a quotient
of the tensor algebra $\T(\gog)$:
\[j^*j_!(\bone_\gh)=\U(\gh^{(1)})\otimes_{\U(\gh)}\bone_\gh=\T(\gog)/\left(J+\T(\gog)\gh\right).\]
Here $J$ denotes the two sided ideal generated by $hg-gh-[h,g]_\gog$ for
$h\in\gh$ and $g\in\gog$. 

\paragraph
There are two natural increasing filtrations on the $\gh$-module
$j^*j_!(\bone_\gh)$.  The first one is induced from the natural
filtration on $\U(\gh^{(1)})$, for which elements of $\gh^{(1)}$ have
degree $1$.  The second one is induced by the grading on $\T(\gog)$,
where elements of $\T^k(\gog)$ have degree $k$.  Throughout this paper
we shall only work with the latter filtration, which shall be denoted
by $F^0\subset F^1\subset F^2 \cdots \subset F^k \cdots$.  Explicitly,
$F^k$ consists of those elements of $j^*j_!(\bone_\gh)$ that have a
lift to $\T(\gog)$ of degree $\leq k$.

\begin{Lemma}
\label{lem:2} 
The associated graded $\gh$-module $\gr(F^\cdot)$ of the above filtration is
precisely $\T(\gn)$.  In other words the successive quotients
$F^k/F^{k-1}$ are isomorphic, as $\gh$-modules, to $\gn^{\otimes k} $.
\end{Lemma}

\begin{Proof}
As $j^*j_!\bone_\gh$ is a quotient of $\T(\gog)$ by the sum of two
ideals, we will understand this quotient in two steps corresponding to
the two succesive quotients. 

We first give a description of the algebra $A = \T(\gog)/J$.  The
ideal $J$ is generated by the linear subspace $R \subset
\gog^{\otimes2}\oplus\gog$ spanned by elements of the form
$hg-gh-[h,g]_{\gog}$ for $h\in\gh,g\in\gog$.  Let $qR$ be the image of
$R$ through the projection
$\gog^{\otimes2}\oplus\gog\to\gog^{\otimes2}$, and form the graded
algebra $qA = \T(\gog)/\langle qR\rangle$, where $\langle qR\rangle$
denotes the two-sided ideal generated by $qR$.  Since $qR$ lies in the
kernel of the quotient algebra morphism $\T(\gog)\to \gr(A)$, we obtain
a surjective algebra morphism $qA\to \gr(A)$.

The quadratic algebra $qA$ is Koszul (see
Appendix~\ref{sec:appendix}), and for such algebras we can apply a
simple criterion~\cite{BG} to check that the map $qA\ra\gr(A)$ is an
isomorphism.  We describe this result below.

Let $\phi:qR\ra \gog$ be the linear map defined as follows. For $x\in
qR$, $\phi(x)$ is the linear part of a preimage of $x$ under the
projection $R\ra qR$.  This is well defined because $R\cap \gog =
0$.  Now Theorem~4.1 in~\cite{BG} states that the morphism $qA\to
\gr(A)$ is an isomorphism if and only if the following conditions are
satisfied:
\begin{itemize}
\item[(1)] $\img(\phi\otimes{\rm id}-{\rm id}\otimes\phi)\subset qR$ 
(this map is defined on $qR\otimes\gog\cap\gog\otimes qR$). 
\item[(2)] $\phi\circ(\phi\otimes{\rm id}-{\rm id}\otimes\phi)=0$. 
\end{itemize}

In our situation, $qR$ is the vector subspace of $\gog\otimes \gog$
spanned by $\{hg-gh|h\in\gh,g\in\gog\}$.  The map $\phi$ maps $hg-gh$
to $[h,g]_{\gog}$.  Thus condition (2) follows from 
the Jacobi identity, while condition (1) is ensured by the stability
of $\gh$ under the bracket of $\gog$.   We conclude that the map $qA
\ra \gr(A)$ is an isomorphism.

In particular, the $\bbk$-vector space $A$ can now be identified with
$\T(\gn)\otimes \S(\gh)$.  Choose $\bbk$-linear splittings of the
projections $\gog\twoheadrightarrow\gn$ and
$\T(\gh)\twoheadrightarrow\S(\gh)$.  Then the composed map
\begin{equation}\label{:composed}
  \T(\gn)\otimes
  \S(\gh)\hookrightarrow\T(\gog)\otimes\T(\gog)\to\T(\gog)\twoheadrightarrow
  A 
\end{equation}
is an isomorphism of filtered $\bbk$-vector spaces. Applying the
counit of $\S(\gh)$ then produces a $\bbk$-linear projection
$\varphi:A\twoheadrightarrow\T(\gn)$.

Let us now prove that $\ker(\varphi)=A\gh$. The kernel of $\varphi$
clearly lies in $A\gh$. Conversely, we now prove that any element
$\sum_sa_sh_s\in A\gh$ lies in the kernel. For any $s$ we can write
$a_s=\sum_t b_tc_t$ in a unique way with $b_t$, resp.~$c_t$, in the
image of $T(\gn)$, resp.~$\S(\gh)$, through \eqref{:composed}. Then
$a_sh_s=\sum_t b_t(c_th_s)\in\ker(\varphi)$ since $c_th_s$ lies in the
augmentation ideal of $\S(\gh)$.

Therefore we get a filtered isomorphism of $\bbk$-vector spaces
$\T(\gn)\tilde\to A/A\gh$ obtained as the composed map
$$
\T(\gn)\hookrightarrow A\twoheadrightarrow A/A\gh\,,
$$
where the first inclusion is determined by a $\bbk$-linear splitting
of $\gog\twoheadrightarrow\gn$.  We now prove that at the level of
associated graded it respects $\gh$-module structures on both sides.
For any $h\in\gh$ and any monomial $x_1\cdots x_k\in\T^k(\gn)$, the
failure of $\gh$-linearity is given by
$$
\sum_{i=1}^kx_1\cdots[h,x_i]_{|\gh}\cdots x_k
$$
where $[\,,]_{|\gh}$ is the $\gh$-part of the bracket, which is
defined by means of the above splitting.  We conclude with the very
simple observation that for any $h\in\gh$ and any
$x_1,\dots,x_k\in\gog$ we have, in $F^{k+1}$,
$$
x_1\cdots h\cdots x_k\in A\gh+F^{k}\,.
$$
Therefore the failure of $\gh$-linearity vanishes after passing to the
associated graded $\gh$-module of $A/A\gh$. 
\footnote{The very same argument shows that the
  $\gh$-module isomorphism $\T^k(\gn)\to F^k/F^{k-1}$ constructed this
  way does not depend on the choice of a $\bbk$-linear splitting
  $\gog\twoheadrightarrow\gn$. }. \qed
\end{Proof}

\paragraph
Next we relate the extension class $\alpha_\gn$ with the filtration $F^\cdot$ on $j^*j_!\bone_\gh$. The inclusion $F^0\hookrightarrow F^k$ of the ground field always
splits for any $k>0$. We shall denote the reduced filtration by
$\tilde{F}^\cdot$.

By the above lemma we have $\tilde{F}^1 \cong \gn$ and
$\tilde{F}^2/\tilde{F}^1 \cong \gn^{\otimes 2}$. Hence the inclusion $\tilde{F}^1\hookrightarrow \tilde{F}^2$ defines a short exact sequence
of $\gh$-modules
\begin{equation}
\label{ses:4} 
0\rightarrow \gn \rightarrow \tilde{F}^2 \rightarrow
\gn^{\otimes 2} \rightarrow 0 .
\end{equation} 
The next lemma shows that the extension class of this sequence is
precisely the class $\alpha := \alpha_\gn\in \Ext^1_\gh(\gn\otimes
\gn, \gn)$ defined in~(\ref{subsec:alpha}).

\begin{Lemma}
\label{lem:3} 
The short exact sequences \eqref{ses:3} and \eqref{ses:4} are
isomorphic and hence both define the same obstruction class $\alpha$.
\end{Lemma}

\begin{Proof} 
We construct a map between 
\[ Q:=\gn\oplus (\gog\otimes \gn) / \langle (h(x),0)-(0,h\otimes
x)\rangle \] 
and $\tilde{F}^2$ which makes all the squares commute. The
required map has two components: one from $\gn$ and the other from
$\gog\otimes \gn$.  The first component is the natural
inclusion map $\gn = \tilde F^1 \subset \tilde F^2$. 
The second one is given by 
\[g\otimes x \mapsto [g\otimes \bar{x}] \] 
where we first choose a lift $\bar{x}$ of $x\in \gn$ to $\gog$ and
then take the class of $g\otimes \bar{x} \in \T^2\gog$ in $\tilde{F}^2$.  A direct computation checkes that the map is well-defined (independent of lifting) and
that the resulting map factors through $Q$. A quick diagram chasing
shows that all squares commute.  \qed
\end{Proof}

\paragraph
The perhaps surprising result that will be proved in
Proposition~\ref{prop:fon} below is that the vanishing of the
extension class $\alpha$, which by the above lemma is only equivalent
to the splitting of the first nontrivial inclusion $F^1\hookrightarrow
F^2$, is in fact equivalent to the splitting of the entire filtration
$F^\cdot$. We will need the following standard lemma which establishes
an isomorphism of $\gog$-modules analogous to the projection formula
in algebraic geometry.

\begin{Lemma}
\label{lem:proj}
Let $i: \gh\hookrightarrow \gog$ be an inclusion of Lie algebras. Let
$E$ be an $\gh$-module and $F$ be a $\gog$-module. Then we have an
isomorphism of $\gog$-modules
\[ i_!(E)\otimes F \cong i_!(E\otimes i^* F).\]
\end{Lemma}
\vspace{-6mm}

\begin{Proof}
Since the result is well-known to the experts, we only provide a short
outline of its proof.  Let $\Delta:\U(\gog)\rightarrow \U(\gog)\otimes
\U(\gog)$ be the cocommutative coproduct on the universal enveloping
algebra, and let $S:\U(\gog)\rightarrow \U(\gog)^\op$ be the antipode
map. We shall freely use the sumless Sweedler notation for the coproduct,
\[ \Delta(f)=f^{(1)}\otimes f^{(2)},\  (\Delta\otimes{\rm
  id})\circ\Delta(f)=({\rm
  id}\otimes\Delta)\circ\Delta(f)=f^{(1)}\otimes f^{(2)}\otimes
f^{(3)}, \  \ldots \]
It is then straightforward to check that the linear map
\[ \phi: i_!(E)\otimes
F=(\U(\gog) \otimes_{\U(\gh)} E)\otimes F \rightarrow i_!(E\otimes i^*
F)=\U(\gog) \otimes_{\U(\gh)} (E\otimes F) \]
given by
\[ \phi((f\otimes x) \otimes y)= f^{(1)}\otimes (x\otimes
S(f^{(2)})y)\] 
is a well-defined isomorphism, with inverse 
\[ \psi: i_!(E\otimes i^*
F)=\U(\gog) \otimes_{\U(\gh)} (E\otimes F) \ra i_!(E)\otimes
F=(\U(\gog) \otimes_{\U(\gh)} E)\otimes F \] 
given by
\begin{equation}
\psi(f \otimes(x\otimes y))= (f^{(1)}\otimes x)\otimes
f^{(2)}y. \tag*{\qed}
\end{equation}
\end{Proof}
\vspace{-6mm}

\begin{Proposition} 
\label{prop:fon}
The following two statements are equivalent:
\begin{itemize}
\item[(a)] The filtration $F^\cdot$ splits.
\item[(b)] The extension class $\alpha$ is trivial.
\end{itemize}
\end{Proposition}

\begin{Proof} 
The implication from (a) to (b) follows from Lemma~\ref{lem:3}. For the other implication we would like to show that the short exact sequences
\[ 0\rightarrow F^{k-1} \rightarrow F^{k} \rightarrow F^{k}/F^{k-1}=\gn^{\otimes k}\rightarrow 0\]
split assuming that the extension class $\alpha$ vanishes. Note that the last equality in the above sequences is proved in Lemma~\ref{lem:2}. Below we will explicitly construct $\gh$-linear maps $\gn^{\otimes k}\rightarrow F^k$ that split the above short exact sequences. 

By Lemma~\ref{lem:obs} the condition $\alpha=0$ is equivalent to the existence of a $\gh^{(1)}$-module structure on $\gn$ that extends the $\gh$-module structure on it. Denote by $\overline{\gn}$ such an extension. Note that as a vector space $\overline{\gn}$ is the same as $\gn$. Denote by $\ad$ the structure map $\gh^{(1)}\otimes \overline{\gn} \rightarrow \overline{\gn}$ for the $\gh^{(1)}$-module $\overline{\gn}$. 

We have a natural map of $\gh$-modules $\gn\hookrightarrow j^*j_! (\bone_\gh)$ for the inclusion of Lie algebras $j:\gh\rightarrow \gh^{(1)}$. 
By adjunction this defines a map of $\gh^{(1)}$-modules
\[ j_!(\gn) \rightarrow j_!(\bone_\gh).\]
 Tensoring both sides with $\overline{\gn}$ yields a map $j_!(\gn)\otimes \overline{\gn} \rightarrow j_!(\bone_\gh)\otimes \overline{\gn}$. 
Applying the projection formula in Lemma~\ref{lem:proj} (for the inclusion $j:\gh\hookrightarrow \gh^{(1)}$) we get a map
\[ j_! (\gn^{\otimes 2}) \rightarrow j_!(n).\]
Iterating this procedure yields for any nonnegative integer $k$ a map of $\gh^{(1)}$-modules
\[ j_! (\gn^{\otimes k+1}) \rightarrow j_!(\gn^{\otimes k}).\]
Hence fixing the integer $k$ we can consider the composition
\[ j_!(\gn^{\otimes k})\rightarrow j_!(\gn^{\otimes k-1}) \rightarrow \cdots \rightarrow j_!(\bone_\gh).\]
Applying adjunction to this composition we get a map of $\gh$-modules
\[ s_k: \gn^{\otimes k} \rightarrow j^*j_!(\bone_\gh).\]
We need to check that the image of $s_k$ lies inside the $k$-th step
of the filtration $F^\cdot$ and that it splits the natural surjective 
map from $F^k$ to $\gn^{\otimes k}$ constructed in Lemma~\ref{lem:2}.

By construction the maps $t_{k+1}: j_!(\gn^{\otimes k+1}) \rightarrow
j_!(\gn^{\otimes k})$ fit into the commutative diagram
\[\begin{CD}
  j_!(\gn^{\otimes k+1}) @>t_{k+1}>> j_!(\gn^{\otimes k}) \\
  @VV\psi V                                  @AA\phi A \\
  j_!(\gn^{\otimes k} )\otimes \overline{\gn} @>t_k\otimes \id>>
  j_!(\gn^{\otimes k-1})\otimes \overline{\gn},\end{CD}\] 
where $\psi$ and $\phi$ are the maps defined in the proof of the
projection formula Lemma~\ref{lem:proj}. This inductive construction
begins with the map $t_1:j_! (\gn) \rightarrow j_!(\bone_\gh)$ which
is explicitly given by $f\otimes x \mapsto fx \otimes 1$. Hence with
respect to the filtrations induced from $\T(\gog)$, $t_1$ increases
the filtration degree by $1$.

It is now important to observe that the coproduct used in the
definition of $\phi$ and $\psi$ not only preserves the filtration for
which elements in $\gh^{(1)}$ are of degree $1$, but also preserves
the filtration induced from $\T(\gog)$. This can be seen by observing
that the natural surjective map $\T(\gog)\rightarrow \U(\gh^{(1)})$ is
a morphism of bialgebras (hence in particular a map of
coalgebras). Thus by an induction argument we conclude that the maps
$t_k$ all increase the filtration degree by $1$. The splitting maps
$s_k$ can then be described as the compositions
\[\begin{CD}
  \gn^{\otimes k} @>\eta_k>> j^*j_!(\gn^{\otimes k}) @>j^*t_k>>
  j^*j_!(\gn^{\otimes k-1}) @>>> \cdots @>>>
  j^*j_!(\bone_\gh).\end{CD}\] 
Here $\eta_k$ is the unit of the adjunction applied to $\gn^{\otimes
  k}$, explicitly given by $ x_1\otimes\cdots\otimes x_k\mapsto
1\otimes_{\U(\gh)} (x_1\otimes\cdots\otimes x_k)$. As $\eta_k$
decreases the filtration by $k$ and we have $k$ times $t$'s to
post-compose with it, the map $s_k$ will end up being a filtration
preserving map, i.e., its image lies inside $F^k\subset
j^*j_!(\bone_\gh)$.

A direct computation of the map $s_k$ shows that it splits the surjective map $F^k\rightarrow \gn^{\otimes k}$. For instance in the cases $k=2$ and $k=3$, we have
\begin{itemize}
\item[k=2]   $x_1\otimes x_2 \mapsto x_1 x_2 - \ad(x_1)x_2.$
\item[k=3]   $x_1\otimes x_2\otimes x_3 \mapsto x_1x_2x_3-x_1\ad(x_2)x_3-x_2\ad(x_1)x_3+\ad(x_2)\ad(x_1)x_3- \ad(x_1)x_2x_3+\ad(\ad(x_1)x_2) x_3.$
\end{itemize}
Lifting is assumed in this formula whenever necessary. One can check
directly that that this formula is independent of all liftings
involved. However it does depend on the choice of the
$\gh^{(1)}$-module $\bar{\gn}$ which lifts the $\gh$-module structure
on $\gn$. \qed
\end{Proof}

%% file: pbw.tex
\section{PBW for inclusions of Lie algebras}
\label{sec:pbw}

In this section we prove the main result of this paper,
Theorem~\ref{thm:mainthm}.  Explicitly, we show that the filtration on
$\U(\gog)/\U(\gog)\gh$ splits if and only if the class $\alpha$
vanishes. 

\paragraph
We shall concentrate our attention on the $\gh$-representation 
\[ i^*i_!(\bone_\gh)=\U(\gog)/\U(\gog)\gh. \] 
This module can be realized as the quotient
$\T(\gog)/\left(I+\T(\gog)\gh\right)$ where $I$ is the two-sided ideal
generated by
\[ \{ g_1\otimes g_2-g_2\otimes g_1-[g_1,g_2]~|~g_1, g_2\in \gog \} \,.\]
We denote by 
\[ R^0\subset R^1\subset\cdots\subset R^k\subset \cdots \] 
the filtration by $\gh$-submodules of $i^*i_!(\bone_\gh)$ induced from
the degree filtration on $\T(\gog)$. We set $G^k:=R^k/R^{k-1}$.

\paragraph
Consider the map
\[ j^*j_!(\bone_\gh)\rightarrow j^*k^*k_!j_!(\bone_\gh)=i^*i_!(\bone_\gh) \]
constructed using the unit map of the adjunction $k_!\dashv
k^*$. This map preserves the filtrations and descends to maps between
associated graded $\gh$-modules
\[ \tau: \T(\gn) = \gr(j^*j_! (\bone_\gh)) \rightarrow
\gr(i^*i_!(\bone_\gh)). \] 
From the descriptions of $j^*j_!(\bone_\gh)$ and $i^*i_!(\bone_\gh)$ via quotients of $\T(\gog)$ we see that the map $\tau$ is surjective.

\begin{Lemma}
\label{lem:4} The kernel of the map $\tau$ is generated by the 
commutators $x\otimes y -y\otimes x$ for $x,y\in\gn$.
\end{Lemma}

\begin{Proof} The idea is to use the classical $\PBW$ theorem for a single Lie algebra. 
Consider the increasing filtration $E^0\subset\cdots\subset 
E^k\subset\cdots$ on the universal enveloping algebra $\U(\gog)$. The
classical $\PBW$ theorem asserts that the kernel of the canonical
surjective map $\gog^{\otimes k}\rightarrow E^k/E^{k-1}$ is generated by the
commutators of elements in $\gog$, thus yielding an isomorphism
between the $k$-th symmetric tensors on $\gog$ and $E^k/E^{k-1}$.

As all these filtrations are compatible (they all arise from the degree
filtration on $\T(\gog)$), the surjective map $j^*j_!(\bone_\gh)\rightarrow i^*i_!(\bone_\gh)$ induces surjections 
on the associated graded to give maps $\gn^{\otimes k} \rightarrow G^k$. Consider the following commutative diagram
\[\begin{CD} 0 @>>> I_1 @>>> \gog^{\otimes k} @>>> E^k/E^{k-1} @>>> 0
\\ @.  @VVV @VVV @VVV @. \\ 0 @>>> I_2 @>>> \gn^{\otimes k} @>\tau^k>> G^k
@>>> 0
\end{CD}\] 
where $I_1$ is the degree $k$ part of the commutator ideal in
$\T(\gog)$ by the $\PBW$ theorem and $I_2$ is the kernel of the
map $\gn^{\otimes k} \rightarrow G^k$.

We want to show that $I_2$ is the $k$-th commutator in $\gn$. It
suffices to show that the map $I_1\rightarrow I_2$ is surjective. By
the snake lemma this is equivalent to showing that the map from the
kernel of $\gog^{\otimes k}\rightarrow \gn^{\otimes k}$ to the kernel
of $E^k/E^{k-1}\rightarrow G^k$ is surjective.

For that we consider another commutative diagram:
\[\begin{CD} 0 @>>> E^{k-1} @>>> E^k @>>> E^k/E^{k-1} @>>> 0 \\ @.
@VVV @VVV @VVV @. \\ 0 @>>> R^{k-1} @>>> R^k @>>> G^k @>>> 0.
\end{CD}\] 
Since all the vertical maps are surjective, the snake lemma shows that we have a surjection
from the kernel of $E^k\rightarrow R^k$ to the kernel of $E^k/E^{k-1}\rightarrow G^k$. But the kernel of $E^k\rightarrow R^k$ is the right ideal
generated by $\gh$ in $\U(\gog)$ intersected with $E^k$, which is a subset of the kernel of $\gog^{\otimes k}\rightarrow \gn^{\otimes k}$. Thus the kernel of $\gog^{\otimes k}\rightarrow \gn^{\otimes k}$ also surjects onto the kernel of $E^k/E^{k-1}\rightarrow G^k$. Thus the lemma is proved.  \qed
\end{Proof}
\medskip

\noindent
To state an if and only if condition for the $\PBW$ property of inclusions of Lie algebras, we need the following lemma concerning the obstruction class $\alpha$.

\begin{Lemma}
\label{lem:5} The obstruction class $\alpha\in \Ext^1(\gn\otimes \gn,
\gn)$ factors through $\S^2(\gn)$.
\end{Lemma}

\begin{Proof} 
The lemma can be seen as a corollary of Lemma~\ref{lem:3} and
Lemma~\ref{lem:4}. Indeed, by Lemma~\ref{lem:3}, we can consider the
following commutative diagram:
\[\begin{CD} 0@>>> \gn @>>> \tilde{F}^2 @>>> \gn\otimes \gn @>>> 0\\
@.  @VVV @VVV @VVV @. \\ 0@>>> \gn @>>> \tilde{R}^2 @>>> G^2 @>>> 0
\end{CD}\] 
where the vertical maps are all defined via the adjunction
$j^*j_!(\bone_\gh)\rightarrow i^*i_!(\bone_\gh)$. By Lemma~\ref{lem:4}
$G^2 = \S^2(\gn)$ and the last vertical map is the canonical quotient
from the tensor product to the symmetric product. Direct calculation shows that the second square is Cartesian. Thus
the lemma is proved.  \qed
\end{Proof}
\medskip

\noindent
We can summarize our main result in the following theorem.

\begin{Theorem} 
Let $\bbk$ be a field and let $\gh\hookrightarrow \gog$ be an inclusion
of Lie algebras over $\bbk$. Consider the two filtrations $R^0\subset
R^1\subset\cdots\subset R^k\subset \cdots $ and $F^0\subset
F^1\subset\cdots\subset F^k\subset \cdots $ defined above. We have:
\begin{itemize}
\item[--] $\gr(F^\cdot) = \T(\gn)$;
\item[--] $\gr(R^\cdot) = \S(\gn)$.
\end{itemize} 
Moreover, if the field $\bbk$ has characteristic zero, then 
the following are equivalent:
\begin{itemize}
\item[(a)] The extension class $\alpha$ is trivial.
\item[(b)] The filtration $F^0\subset F^1\subset\cdots\subset
F^k\subset \cdots $ splits;
\item[(c)] The filtration $R^0\subset R^1\subset\cdots\subset
R^k\subset \cdots $ splits.
\end{itemize} 
In fact, if the extension class $\alpha$ is trivial, we have the
following explicit splitting of the filtration $R^\cdot$ that resembles the
standard $\PBW$ isomorphism:
\[ I: \S(\gn)\rightarrow \T(\gn)\cong j^*j_!(\bone_\gh)
\rightarrow i^*i_!(\bone_\gh)\cong \U(\gog)/\U(\gog)\gh. \] 
Here the first arrow is given by any (graded) splitting of the surjective morphism $\T(\gn)\to\S(\gn)$ in $\gh$-$\mathfrak{Mod}$. 
\end{Theorem} 

\begin{Proof}
The fact that $\gr(F^\cdot)=\T(\gn)$ is proved in Lemma~\ref{lem:2}, and $\gr(R^\cdot)=\S(\gn)$ follows from Lemma~\ref{lem:4}. For the second part of the theorem, 
Proposition~\ref{prop:fon} shows that (a) and (b) are equivalent. By Lemma~\ref{lem:5} (c) implies (a). Below we will show that (b) implies (c) and hence 
all of (a), (b), (c) are equivalent.

Assuming a splitting $I$ of the surjection $\T(\gn)\rightarrow
\S(\gn)$ (which always exists over a field of characteristic zero) and a splitting $s$ of the filtration $F^\cdot$, 
we can define the following composition
\[ \S^k(\gn) \stackrel{I^k}{\rightarrow} \T^k(\gn) \stackrel{s^k}{\rightarrow} F^k \rightarrow R^k.\]
Here the last map is the canonical surjective map from $F^k$ to $R^k$. The fact that this composition defines a splitting for the filtration $R^\cdot$ follows 
from the following commutative diagram
\[\begin{CD} \S^k(\gn) @>I^k>> \T^k(\gn) = F^k/F^{k-1}
@>s^k>> F^k\\ @.  @VVV
@VVV \\ @.  \S^k(\gn) = R^k/R^{k-1} @<<< R^k.
\end{CD}\] 
\vspace{-2.37em}

\qed
\end{Proof}

%% file: twisted.tex
\section{Generalization to any representation}
\label{sec:twisted}

The main goal of this section is to extend Theorem~\ref{thm:mainthm}
from the case of the trivial representation ${\bf 1}_\gh$ to that of
an arbitrary finite dimensional $\gh$-representation $V$.  Consider
the filtrations $R^0\subset R^1\subset\cdots R^k\subset\cdots$,
$F^0\subset F^1\subset\cdots F^k\subset\cdots$, on $i^*i_!(V)$,
$j^*j_!(V)$, respectively, which are induced by the degree filtration
on $T(\gog)\otimes V$.  Then we have the following theorem.

\begin{Theorem}
\label{thm:twisted}
There are isomorphisms of $\gh$-modules $gr(F^*)=T(\gn)\otimes V$ and
$gr(R^*)=S(\gn)\otimes V$.  Moreover, the following are equivalent:
\begin{itemize}
\item[(a)] The extension classes $\alpha$ and $\alpha_V$ are trivial. 
\item[(b)] The filtration $F^0\subset F^1\subset\cdots
  F^k\subset\cdots$ splits.
\item[(c)] The filtration $R^0\subset R^1\subset\cdots
  R^k\subset\cdots$ splits.
\end{itemize}
\end{Theorem}

\noindent\emph{Sketch of proof.} 
First observe that as vector spaces we have filtered $\bbk$-linear
isomorphisms
$$
i^*i_!(V)\cong i^*i_!({\bf 1}_\gh)\otimes V\quad\textrm{and}\quad j^*j_!(V)\cong j^*j_!({\bf 1}_\gh)\otimes V
$$
They are not isomorphisms of $\gh$-modules, but on the level of associated graded they induce $\gh$-module isomorphisms. 
This proves the first part of the theorem.

Contrary to the situation when the representation is trivial, the
inclusions $F^0\hookrightarrow F^k$ and $R^0\hookrightarrow R^k$ do
not automatically split. In particular the inclusion
$V=F^0=R^0\hookrightarrow F^1=R^1$ splits if and only if $\alpha_V$ is
trivial.

Finally, if $\alpha_V$ is trivial then there exists an $\gh^{(1)}$-module $\tilde{V}$ such that 
$Res(\tilde{V})=V$. From this we deduce an isomorphism of $\gh$-modules 
$j^*j_!(V)\cong j^*j_!({\bf 1}_\gh)\otimes V$. We conclude by using the fact that the Theorem is true 
for the trivial representation. \qed

%% file: example.tex
\section{An example of a non trivial class}
\label{ex}

We now give an example of an inclusion of Lie algebras
$\mathfrak{h}\hookrightarrow\mathfrak{g}$ for which the obstruction
class is non trivial.  Let $\mathfrak{g}=\mathfrak{sl}_2$; recall that
it is generated by $e$, $h$ and $f$, satisfying the relations
$$
[e,f]=h\,,\quad[h,e]=2e\,,\quad[h,f]=-2f\,. 
$$
Now let $\mathfrak{h}$ be the Lie subalgebra generated by $e$ and $h$. Then $\mathfrak{n}=\mathfrak{g}/\mathfrak{h}$ is the one dimensional $\mathfrak{h}$-module generated as a vector space by $f$, with module structure defined by 
$$
e\cdot f=0\quad\textrm{and}\quad h\cdot f=-2f\,.
$$
\begin{Proposition}
The obstruction class $\alpha$ is non-trivial. 
\end{Proposition}

\begin{Proof}
First observe that the Chevalley-Eilenberg 1-cocycle 
\[ c\in C^1\big(\mathfrak{h},\Hom(\mathfrak n,\mathfrak h)\big) \]
given by 
$$
c(e)(f)=e \cdot f-[e,f]=-h,\quad c(h)(f)=h\cdot f-[h,f]=0
$$
is a representative of the exact sequence 
$$
0\ra \mathfrak h\ra\mathfrak g\ra\mathfrak n\ra 0\,.
$$
Therefore the $1$-cocycle $a\in C^1\big(\mathfrak{h},{\rm Hom}(\mathfrak n^{\otimes2},\mathfrak n)\big)$ given by 
$$
a(e)(f,f)=-h\cdot f=2f,\quad a(h)(f,f)=0
$$
is a representative of the obstruction class $\alpha$. 

Finally, observe that since $e$ acts trivially on $\mathfrak n$, then it acts trivially on ${\rm Hom}(\mathfrak n^{\otimes2},\mathfrak n)$. Consequently, 
for any $b\in{\rm Hom}(\mathfrak n^{\otimes2},\mathfrak n)$ we have
$d(b)(e)=0$, so that $a\neq d(b)$.  Thus $\alpha\neq0$. 
\qed
\end{Proof}

%% file: appendix.tex
\appendix

\section{The algebra {\em qA} is Koszul}
\label{sec:appendix}

In this appendix we prove that the quadratic algebra $qA$, defined as
the quotient of $\T(\gog)$ by the two sided ideal generated by the
linear subspace $qR$ of $\gog^{\otimes 2}$ spanned by
$\{hg-gh~|~h\in\gh,g\in\gog\}$, is Koszul.  We refer to \cite{BG} and
references therein for the many definitions of Koszulity and their
equivalence.

The Koszul complex $K(qA)$ of $qA$ is a subcomplex of the Bar resolution $B_{qA}(\bbk)$ of $\bbk$ 
as a left $qA$-module via the augmentation map $\varepsilon:qA\to\bbk$:  
$$
K(qA):=\bigoplus_{i\geq0}\left(qA\otimes \tilde{K}^i(qA)\right)[i]
\subset\bigoplus_{i\geq0}\left(qA\otimes qA^{\otimes i}\right)[i]=: B_{qA}(\bbk),
$$
where 
$$
\tilde{K}^i(qA):=\bigcap_{k=0}^{i-2}\gog^{\otimes k}\otimes qR\otimes\gog^{\otimes i-k-2}.
$$
Recall that the differential on the Bar resolution $B_{qA}(\bbk)$ is defined by 
\[ d(a_0\otimes\cdots\otimes a_i)=\sum_{k=0}^{i-1}(-1)^ka_0\otimes\cdots\otimes a_ka_{k+1}\otimes\cdots\otimes a_i+(-1)^ia_0\otimes\cdots a_{i-1}\varepsilon(a_i).
\]
\begin{Proposition}
For any $i<0$, $H^i\big(K(qA)\big)=0$. In other words the algebra $qA$ is Koszul. 
\end{Proposition}

\begin{Proof}
For any $i>0$, $\tilde{K}^i(qA)=\big(\wedge^{i-1}\gh\big)\wedge\gog$ is the image of 
$\gh^{\otimes(i-1)}\otimes\gog$ in $\gog^{\otimes i}\subset qA^{\otimes i}$ through the total antisymmetrization map 
$$
x_1\otimes\cdots\otimes x_i\mapsto 
x_1\wedge\dots\wedge x_i:=\sum_{\sigma\in S_i}\epsilon(\sigma)x_{\sigma(1)}\otimes\cdots\otimes x_{\sigma(i)}.
$$
As usual, $\tilde K^0(qA)=\bbk$. 
Now observe that the only non-zero term in the restriction of the
differential of the Bar resolution to 
$K(qA)$ is the first one: 
$$
d(\sum_s a_0^{(s)}\otimes\cdots\otimes a_i^{(s)})=
\sum_s
a_0^{(s)}a_1^{(s)}\otimes\cdots\otimes a_i^{(s)}.
$$
This is a general fact that is not specific to the peculiar situation
we are working at.  

In other words, for an element $a\otimes x_1\wedge\dots\wedge x_i\in
qA\otimes\tilde{K}^i(qA)$ we have
\[ 
d(a\otimes x_1\wedge\dots\wedge x_i)
=\sum_{j=1}^i(-1)^{j-1}ax_j\otimes x_1\wedge\dots\wedge\hat{x_j}\wedge\dots\wedge x_i.
\]
Remember that the symmetric algebra $S(\gh)$ is Koszul: its Koszul complex $K\big(S(\gh)\big)$, which is $\oplus_{i\geq0}S(\gh)\otimes\wedge^i(\gh)[i]$ 
with differential being given by the formula above, is acyclic in negative degrees. 
Finally, one sees that $K(qA)$ is isomorphic to the $dg$ $K\big(S(\gh)\big)$-module freely generated 
by the two step complex
$$
\cdots\longrightarrow0\longrightarrow\big(\T(\gn)\otimes\gn\big)[1]\longrightarrow\T(\gn)\longrightarrow0\longrightarrow\cdots,
$$ 
which is acyclic in negative degrees.  This can be easily seen by
considering a $\bbk$-linear splitting $\gn\hookrightarrow\gog$ 
of $\gog\twoheadrightarrow\gn$ and observing that 
\begin{align*}
K\big(S(\gh)\big)\otimes\big(\T(\gn)\otimes(\gn[1]\oplus\bbk)\big)
& \cong
qA\otimes\left(\bigoplus_{i\geq0}\wedge^i(\gh)[i]\right)\otimes(\gn[1]\oplus\bbk)
\\
& \cong qA\otimes\left(\bbk\oplus\bigoplus_{i\geq0}\big(\wedge^{i-1}(\gh)\big)\wedge\gog[i]\right).
\end{align*}
We leave to the reader the straightforward task of tracking the
differential through this identification.  \qed
\end{Proof}